%Paper: alg-geom/9502008
%From: OTTAVIANI@aquila.infn.it
%Date: Mon, 13 Feb 1995 18:18:43 +0100 (WET)

\magnification =\magstep1
\baselineskip =13pt
\def\pp{I\!\! P}
\centerline{\bf  ON  SINGULARITIES OF  ${\cal M}_{I\!\! P^3}(c_1,c_2)$}
 \vskip 1cm

\centerline {\bf Vincenzo  Ancona  and Giorgio Ottaviani  \footnote{*}{\rm
Both authors were  supported by  MURST and by GNSAGA of CNR}}
\vskip 0.8cm

\centerline{alg-geom/9502008}
\vskip 1cm
{\bf  Summary}. Let ${\cal M}_{I\!\! P^3}(c_1,c_2)$ be  the moduli space
of stable rank-$2$ vector bundles on $I\!\! P^3$ with Chern classes $c_1$, $
c_2$.
We prove the following results.

\noindent 1) Let  $0 \le \beta < \gamma $  be two integers, ($\gamma \ge 2)$,
such that
 $2\gamma-3\beta>0$; then ${\cal M}_{I\!\! P^3}(0,2\gamma^2-3\beta^2)$ is
singular
(the case $\beta =0$ was previously proved by M. Maggesi).

\noindent 2) Let $0 \le \beta < \gamma $  be two odd integers ($\gamma \ge 5)$,
 such that $2\gamma-3\beta+1>0$; then ${\cal M}_{I\!\!
P^3}(-1,2(\gamma/2)^2-3(\beta/2)^2+1/4)$
 is singular.

\noindent In particular ${\cal M}_{I\!\! P^3}(0,5)$, ${\cal M}_{I\!\!
P^3}(-1,6)$ are singular.
\vskip 1.5cm

\def\cc{\ I\!\!\!\!C}
\def\pp{I\!\! P}

The first examples of singular moduli spaces of stable vector bundles on a
projective space
 were found by the authors in [AO1], where it is shown that the symplectic
special instanton
 bundles on $\pp^5$ with second Chern class $c_2=3,4$ correspond to singular
points of their moduli space.
 Later R. M. Mir\'o-Roig  [MR]  detected an example in the case of rank-$3$
vector bundle
 on $\pp^3$.
Recently M. Maggesi [M]  has pointed out that pulling back some particular
instanton
bundles with second Chern class $c_2=2$ by a finite morphism $ \pp^3
\rightarrow \pp^3$
 one always obtains singular points in the corresponding moduli spaces, thus
giving examples
in the rank-$2$ case. His result is that for any integer $d\ge 2$ the moduli
space
${\cal M}_{\pp^3}(0,2d^2)$ is singular; $d=2$ gives ${\cal M}_{\pp^3}(0,8)$.

The aim of this paper is to exhibit some more examples. The idea is to replace
the
usual pull-back by the more general construction of "pulling back over
$\cc^4\backslash 0$ ",
 introduced in [Ho] and  developed in [AO2]. The result is the following:
 let  $0 \le \beta < \gamma $  be two integers, ($\gamma \ge 2)$, such that
$2\gamma-3\beta>0$; then
the moduli space ${\cal M}_{\pp^3}(0,2\gamma^2-3\beta^2)$ is singular (Main
theorem I). In particular taking  $\beta = 1,\gamma =2$ we obtain that ${\cal
M}_{\pp^3}(0,5)$
 is singular; in this  case the singular points we have detected fall in the
closure of
the open set consisting of instanton bundles [Ra].

We prove a similar result in the case $c_1=-1$ (Main Theorem II). In particular
we find that
 ${\cal M}_{\pp^3}(-1,6)$ is singular.
\vskip.3cm
\centerline { \bf  \S  1}
\vskip.3cm
Let $U$ be a $2$-dimensional complex vector space, $V=U\oplus U$, and
 $\pp^3=\pp(U\oplus U)$ the projective space of hyperplanes in $V$.
We denote by ${\cal M}_{\pp^3}(c_1,c_2)$ the moduli space of stable rank-$2$
vector bundles
 on $\pp^3$
with Chern classes $c_1$, $ c_2$.

 Let us  choose homogeneous
 coordinates $(a,b,c,d)$ in $\pp^3$ so that $(a,c)$ and $(b,d)$ are coordinates
in  $U\oplus0$ and
 $0\oplus U$ respectively.

It is well known ([H],[LP],[NT]) that any stable rank-$2$ vector bundle $E$  on
$\pp^3$
 belonging to ${\cal M}_{\pp^3}(0,2)$
 can be described by a monad

\vskip.2cm

$$ 0\rightarrow  {\cc^2 \otimes \cal O}_{\pp^3} (-1) \buildrel A \over
 \rightarrow U^3\otimes{\cal O}_{\pp^3}
\buildrel {B}\over\rightarrow \ \cc^2 \otimes {\cal O}_{\pp^3}(1) \rightarrow 0
\eqno (1.1)$$

\vskip.2cm

\noindent where  $$A= \left(\matrix{a&b&0&c&d&0\cr 0&a&b&0&c&d}\right)$$

\vskip.2cm

$$B= \left( \matrix { & &-d&-c\cr &-d&-c& &\cr-d&-c& & \cr & &b&a&\cr &b&a& \cr
b&a& & }\right)
\left( \matrix
{\alpha_3&\alpha_4\cr\alpha_2&\alpha_3\cr\alpha_1&\alpha_2\cr\alpha_0&\alpha_1}
\right)$$

\vskip.2cm

\noindent  and $\alpha_0\dots \alpha_4$ are constant coefficients subject to
the condition

\vskip.2cm

$$det {\left( \matrix
{\alpha_0&\alpha_1&\alpha_2\cr\alpha_1&\alpha_2&\alpha_3\cr\alpha_2
&\alpha_3&\alpha_4}\right)}\ne 0 $$
\vskip.2cm

There is a natural action of $SL(2) \cong  SL(U) $ on the matrices $A$ and $B$
through the transformations

\vskip.2cm

$$\left(\matrix {a \cr c}\right)\rightarrow g \left(\matrix {a \cr c}\right) ,
 \left(\matrix {b \cr d}\right)\rightarrow g \left(\matrix {b \cr d}\right) \ (
g\in SL(2)) $$
\vskip.2cm
It follows that  $SL(2) \cong  SL(U) $ acts on the monad $(1.1) $, hence on its
cohomology $E$.
 More precisely, for $g= \left(\matrix {x&y \cr z&w}\right) $ let $$Q_g=\left(
\matrix { x& & &z& &
\cr  &x& & &z& \cr & & x& & &z \cr y& & &w& &
\cr  &y& & &w& \cr & & y& & &w}\right)$$
\vskip.2cm
It is easy to check that $g^*A= AQ_g$, $g^*B=Q_g^{-1}B$, which implies that the
monads
$(A,B)  $  and $ ( g^*A, g^*B)$ are equivalent. Moreover, the monad $(1.1)$ is
$SL(U)$-invariant
($ SL(U)$ acts trivially on  $\cc^2 $ and diagonally on $U^3$).

\proclaim Lemma 1.
The bundle $E$ defined by the monad $(1.1)$ admits the following
 $SL(U)$-invariant minimal resolution:
 $$ { 0\rightarrow {U\otimes\cal O}_{\pp^3}(-4)\rightarrow (S^2U)^2\otimes
{\cal O}_{\pp^3} (-3)}$$ $$
{\rightarrow [S^3U \otimes {\cal O}_{\pp^3}(-2)}]\oplus
   [\cc^2 \otimes{\cal O}_{\pp^3}(-1)]
 \rightarrow  E \rightarrow 0 \eqno (1.2)$$
\vskip.2cm

{\sl Proof:} The minimal resolution of $E$ is well known and is clearly
$SL(U)$-invariant;
(1.2)  follows by inspecting in it the cohomology groups of $E$ as
$SL(U)$-representations.
\vskip.2cm

\proclaim Lemma 2.
$ H^0E(1)= \cc^2$, $H^0E(2)=S^3(U)\oplus U^2$
as $SL(U)$-representations;
 moreover, for $t\ge 3$    $H^0E(t)$ is a direct sum of symmetric powers
$S^m(U)$ with
$0\le m \le {t+1}$.

\proclaim  Lemma 3.
 $H^1End\ E(-1)=U^4$, $ H^1End\ E= \cc^4\oplus (S^2U)^3$, $ H^1End\
E(1)=U^2\oplus (S^3U)^2$,
 $H^1End\ E(2)=S^4U$  as $SL(U)$-representations.
Moreover,  $\ H^1End\ E(t)=0 $ for $t\le -3$ and $t\ge 3$, and $H^1End\ E(-2)$
is equal to
 $\cc $ if $I=\alpha_0\alpha_4-4\alpha_1\alpha_3+3\alpha_2^2=0$ , to $0$
otherwise.

\vskip.2cm

{\sl  Proof.} The statement about $H^1End\ E(-2)$ can be found in [LP],[NT].
The remaining equalities
 can be easily obtained by tensoring the sequence $(1.2)$ by $E(t)$ and
inspecting
the corresponding cohomology sequences.
\vskip.3cm

\centerline {\bf \S 2}

\vskip.3cm

Let $0 \le \beta < \gamma $  be two integers, ($\gamma \ge 2)$; let  $f_1,\dots
,f_4$
be homogeneous polynomials
 in the variables $a, b, c, d$ without common zeroes of degree $ \gamma
-\beta,\gamma -\beta,
\gamma +\beta,\gamma +\beta$ respectively. Let us take into account the diagram

$$\matrix {\ \cc^4 \backslash 0 &\buildrel {\omega}\over \rightarrow&\  \cc^4
\backslash 0 \cr \eta \downarrow & &\eta \downarrow  \cr \ \ \ \pp^3& &\ \ \
\pp^3  }$$
where $\omega$ is defined by  $f_1,\dots ,f_4$.

According to [Ho], [AO2], from any bundle $E$ defined by the monad $(1.1)$ we
construct
a rank-$2$ bundle $E_{\beta, \gamma } $ such that  $\eta^* E_{\beta, \gamma
}=\omega^*\eta^*E$ . The bundle
$E_{\beta, \gamma } $ is the cohomology of the monad
\vskip.2cm
$$ 0\rightarrow  {\ \cc^2 \otimes {\cal O}_{\pp^3} (-\gamma)
\buildrel A_{\beta, \gamma } \over
 \rightarrow {\cal U}^3
\buildrel B_{\beta, \gamma }\over\rightarrow {\ \cc^2 \otimes {\cal O}}_{
\pp^3}
(\gamma)
  \rightarrow 0}  \eqno (2.1)$$
\vskip.2cm

\noindent where  ${ \cal U}={\cal O}_{\pp^3}(-\beta)
\oplus {\cal O}_{ \pp^3}(\beta)$, and $ A_{\beta, \gamma }$, $B_{\beta, \gamma
}$
\  are  obtained from the matrices $A$, $B$ in $(1.1)$ replacing
$a,b,c,d$ by $f_1,f_2,f_3,f_4 $ respectively.

\noindent In particular the Chern classes of $E_{\beta, \gamma } $ are
 $c_1=0,c_2=2\gamma^2-3\beta^2$.

\noindent Of course $E_{\beta, \gamma }$ depends on $f_1,\dots ,f_4$ but for
simplicity
we omit this fact in the notations.

The cohomology groups of $E_{\beta, \gamma }(t) $ can be computed by  [AO2
 \S 2 ]; in particular, the theorem 2  of [AO2] can be rephrased as follows.
  \vskip.2cm

  \proclaim Theorem 4. Let $H^iE(t)=T^t(U)$ where $T^t$ is a representation of
$SL(U)$. Then $$ h^iE_{\beta,\gamma}(t)= \sum_{h\in {\bf Z}} \sum_{j=0}^4
(-1)^j h^0[\bigwedge^j
\left( {\cal U}^2(-\gamma)\right) \otimes T^h({\cal U})\otimes {\cal O}_{
\pp^3}
(t-h\gamma)]$$

 In practice the groups $ h^iE_{\beta,\gamma}(t)$ can be computed as follows.
Let $s_p$
the dimension of the  the degree $p$
summand of the artinian algebra  $S=\cc [a,b,c,d]/(f_1,f_2,f_3,f_4)$; for $ h
\in Z$ we  write
 $$ T^h({\cal U})(h\gamma)=\bigoplus_s {\cal O}_{\pp^3}(\mu_{s,h}) \eqno (2.2)
$$
 \noindent Let
$$b_q(E)=\sharp \ \{(s,h):\mu_{s,h}=q \} $$
Then
$$ h^iE_{\beta,\gamma}(t)=\sum_{p+q=t} s_p  b_q(E)  \eqno (2.3) $$

\noindent Moreover, $$s_p= \sum_{j=0}^4 (-1)^j h^0[\bigwedge^j
  {\cal U}^2 \otimes {\cal O}_{ \pp^3}
(p-j\gamma)]  \eqno (2.4)$$

The formulae $(2.3), (2.4)$ are  easy adaptations of the results of  [AO2 \S
2].

The above formulae still hold if we replace $E_{\beta,\gamma}$ (resp $E$)
by $End \ E_{\beta,\gamma}$ (resp $End \ E$).

\proclaim Proposition 5. Let $E$ be defined by the monad $(1.1)$. Then
\item {(i)} $ E_{\beta,\gamma}$ is stable if and only if $2\gamma-3\beta >0$.
\item {(ii)} $h^1End \ E_{\beta,\gamma}=k_{\beta,\gamma}+s_{2\gamma} h^1End \
E(-2)$
 where $k_{\beta,\gamma}$ is a constant not depending on $E$.
\vskip.2cm

{\sl Proof}.
 (i) We need to show that $h^0E_{\beta,\gamma}=0$ iff
$2\gamma-3\beta >0$ . We compute $h^0E_{\beta,\gamma}$ substituting $t=0$  in
the formula
 $(2.3)$ .
 Since $h^0E(h)=0$ for $h\le 0$ and $s_p=0$ for $p<0$, a contribution
to the right-hand side of $(2.3)$ can occur only for $q\le 0$. By  lemma 2 the
integers
 $\mu_{s,h}$ appearing in $(2.2)$ are strictly positive for $h=1$, while for
$h=2$
we have $ T^2({\cal U})(2\gamma)= S^3({\cal U})(2\gamma)\oplus ({\cal U}^2)
(2\gamma)=\oplus_s{\cal O}_{\pp^3}(\mu_{s,2})$ with $\inf_s \{  \mu_{s,2}\}=
2\gamma-3\beta$;
for $h\ge 3$  we have $2\gamma-3\beta < h\gamma-(h+1)\beta \le \inf_s \{
\mu_{s,h}\}$;
 the conclusion  follows.

(ii) Again this follows from the formula $(2.3)$ with  $End \ E$ in place of
$E$ and from
the lemma 3.
\proclaim Main Theorem I. Let $0 \le \beta < \gamma $  be two integers ($\gamma
\ge 2)$,
 such that $2\gamma-3\beta>0$. Let $E$
be a bundle belonging to ${\cal M}_{\pp^3}(0,2)$ such that $h^1End \ E(-2)\ne
0$.
The moduli space ${\cal M}_{\pp^3}(0,2\gamma^2-3\beta^2)$ is singular at the
points corresponding
to the bundles $E_{\beta,\gamma}$.
\vskip.2cm

{\sl Proof}. By proposition 5,(i)  we can define a natural map
\vskip.2cm
$$\matrix {{\cal M}_{\pp^3}(0,2)& \rightarrow &{\cal
M}_{\pp^3}(0,2\gamma^2-3\beta^2)
\cr [F] &\rightarrow & [F_{\beta,\gamma}]}$$
\vskip.2cm
\noindent which is  clearly algebraic.

\noindent Let $E$ be as in the statement. In any neighborhood of its class $[E]
\in {\cal M}_{\pp^3}(0,2)$
 there is a $[F]$ with  $h^1End \ F(-2)=0$. From proposition 5 (ii) we obtain
$$h^1End \ E_{\beta,\gamma}  >h^1End \ F_{\beta,\gamma} $$
which clearly implies that $[E_{\beta,\gamma}]$ is a singular point of
${\cal M}_{\pp^3}(0,2\gamma^2-3\beta^2)$.

Taking $\gamma = 2, \beta  =1$ we obtain:
\proclaim Corollary. ${\cal M}_{\pp^3}(0,5)$ is singular.

In this particular case the singular points we have detected fall in the
closure of the open set
consisting of instanton bundles (see [Ra]). It would be interesting to know
whether
the same property  is true in the general case.

Let us remark that taking $\beta = 0$ we recover the result of [M].

By improving our technique we can show more generally that
${\cal M}_{\pp^3}(0,k\gamma^2-(k+1)\beta^2)$ is singular for $0<\beta<\gamma$
and any integer
 $k>0$ such that $k\gamma-(k+1)\beta>0 $ [AO3].
 \vskip.3cm
\centerline { \bf  \S  3}
\vskip.3cm
In this section we deal with singularities of the moduli spaces  ${\cal
M}_{\pp^3}(-1,c_2)$.
 Let $0 < \beta < \gamma $   be two {\sl odd} integers; let  $f_1,\dots ,f_4$
be homogeneous
 polynomials in the variables $a, b, c, d$ without common zeroes of degree
$ (\gamma -\beta)/2, (\gamma -\beta)/2 ,(\gamma +\beta)/2,(\gamma +\beta)/2$
respectively. We can construct the monad
\vskip.2cm
$$ 0\rightarrow  {\ \cc^2 \otimes {\cal O}_{\pp^3} (-(\gamma+1)/2)
\buildrel A_{\beta/2, \gamma/2 } \over
 \rightarrow {\cal W}^3
\buildrel B_{\beta/2, \gamma/2 }\over\rightarrow {\ \cc^2 \otimes {\cal O}}_{
\pp^3}
((\gamma-1)/2)
  \rightarrow 0}  \eqno (3.1)$$
\vskip.2cm
\noindent where  ${ \cal W}={\cal O}_{ \pp^3}(-(\beta+1)/2)
\oplus {\cal O}_{\pp^3}((\beta-1)/2) $, and $ A_{\beta/2, \gamma/2 }$,
 $ B_{\beta/2, \gamma/2 }$
\  are  obtained from the matrices $A$, $B$ in $(1.1)$ replacing
$a,b,c,d$ by $f_1,f_2,f_3,f_4 $ respectively.

\noindent The cohomology of the monad (3.1) is a rank-$2$ bundle
$E_{\beta/2,\gamma/2 } $  whose  Chern classes  are
 $c_1=-1,c_2=2(\gamma/2)^2-3(\beta/2)^2+1/4 $.

Let $E_{\beta, \gamma}$ be defined  as in \S 2 by the monad (2.1), via the
homogeneous
polynomials $f_1(a^2,b^2,c^2,d^2),\dots ,f_4(a^2,b^2,c^2,d^2) $. Then for
${t\in {\bf Z}}$
 $$\pi^* E_{\beta/2,\gamma/2}(t)= E_{\beta, \gamma} (2t-1) \eqno (3.2)$$
where $\pi : \pp^3 \rightarrow \pp^3$ is a finite morphism of degree $8$.
Moreover
$$\pi^* End \ E_{\beta/2,\gamma/2}(t)= End \ E_{\beta, \gamma}(2t) \eqno
(3.3)$$

In order to compute the cohomology groups of
of $ End \ E_{\beta/2,\gamma/2}(t) $ let  $H^1End \ E(h)=R^h(U)$, where  $R^h$
is a representation of $SL(U)$; let ${\cal V}={\cal O}_{ \pp^3}(-\beta/2)
\oplus {\cal O}_{\pp^3}(\beta/2) $. Though ${\cal V}$ is only a ${\bf
Q}$-bundle, by lemma 3
 one easily checks that $R^h({\cal V})(h\gamma/2)=\bigoplus_s {\cal
O}_{\pp^3}(\nu_{s,h})$
 with $\nu_{s,h}\in {\bf Z}$. As in \S 2 we define  $s_p$ as
the dimension of the  degree $p$ summand of the artinian algebra
$S=\cc [a,b,c,d]/(f_1,f_2,f_3,f_4)$; let
$$b_q(End\ E)=\sharp \ \{(s,h):\nu_{s,h}=q \} $$
Then
\proclaim Lemma 6.
 \item {(i)} $ h^1End\ E_{\beta/2,\gamma/2}(t)=\sum_{p+q=t} s_p  b_q(End\ E) $
\item {(ii)} $h^1End \ E_{\beta/2,\gamma/2}=k_{\beta/2,\gamma/2}+s_{\gamma}
h^1End \  E(-2)$
 where $k_{\beta/2,\gamma/2}$ is a constant not depending on $E$.
\vskip.2cm

{\sl Proof}. (The details are left to the reader). By (3.3)  the groups $
h^1End\ E_{\beta/2,\gamma/2}(t)$, ($t\in {\bf Z}$),
 are uniquely determined by the groups $ h^1End\ E_{\beta,\gamma}(h)$, ($h\in
{\bf Z}$);
by comparison with   the formula (2.3) (applied to $ End \ E_{\beta,\gamma}$)
one checks (i);
 (ii) is an  immediate consequence of (i).

\proclaim Main Theorem II. Let $0 < \beta < \gamma $  be two odd integers
($\gamma \ge 5)$,
 such that $2\gamma-3\beta+1>0$. Let $E$
be a bundle belonging to ${\cal M}_{\pp^3}(0,2)$ such that $h^1End \ E(-2)\ne
0$.
The moduli space ${\cal M}_{\pp^3}(-1,2(\gamma/2)^2-3(\beta/2)^2+1/4) $ is
singular at the point corresponding
to the bundle $E_{\beta/2,\gamma/2}$.

The proof is similar to the proof of the main theorem I. The condition
$2\gamma-3\beta+1>0$
ensures by (3.2) the stability of $E_{\beta/2,\gamma/2}$, while the condition
$\gamma \ge 5$
implies $s_{\gamma}\ne 0$ in lemma 6.

Taking  $\gamma = 5$, $\beta = 3$ we obtain:

 \proclaim Corollary. ${\cal M}_{\pp^3}(-1,6)$ is singular.

\vskip1.cm

\centerline {\bf REFERENCES}

 \item{[AO1]} V. Ancona, G. Ottaviani. On moduli of instanton bundles on
$\pp^{2n+1}$ .
  {\it   Pacific \ Journal \ of  \ math.} (to appear).
 \item{[AO2]} V. Ancona, G. Ottaviani. The Horrocks bundles of rank three on
${\pp^5}$.
{\it J.  reine \ angew.\ Math.}  {\bf 460}  (1995).
\item{[AO3]} V. Ancona, G. Ottaviani. In preparation.
\item{[H]}  R. Hartshorne. Stable vector bundles of rank $2$ on $ \pp^3$. {\it
Math. Ann. }
{\bf  238} (1978), 229-280.
\item{[Ho]} G. Horrocks. Examples of rank three vector bundles on
five-dimensional
 projective space. {\it J. London \ Math. \ Soc. } {\bf  18}  (1978), 15-27.
\item{[LP]} J. Le Potier. Sur l'espace des  modules des fibr\'es de Yang et
Mills. In
Math\'emathique  et Physique, {\it  S\'em.\  Ecole \  Norm. \ Sup. 1979-1982},
65-137,
 Basel-Stuttgart-Boston (1983).
\item{[M]} M. Maggesi.   Tesi \ di \ laurea. {\it Universita' \ di \ Firenze}
(1994).
\item{[MR]} R.M. Mir\'o-Roig. Singular moduli space of stable vector bundles on
$\pp^3$.
 {\it   Pacific \ Journal \ of  \ math.} (to appear).
\item{[NT]} M.S. Narasimhan, G. Trautmann. The Picard group of the
compactification of

  ${\cal M}_{\pp^3}(0,2)$.  {\it  J.  reine \ angew.\ Math.},
{\bf 422} (1991), 21-44.
\item{[Ra]} A.P. Rao. A family of vector bundles on $\pp^3$. In {\it Lecture \
Notes \ in \ Math.}
{\bf 1266} 208-231, Berlin, Heidelberg, New York: Springer 1987.
\vskip 2cm

Authors' addresses:

\noindent V.A.: Dipartimento di Matematica "U. Dini"

\noindent Viale Morgagni 67/A

\noindent I-50134 Firenze (Italy)

\noindent e-mail ancona@udini.math.unifi.it
\vskip .3cm
\noindent G.O.: Dipartimento di Matematica

\noindent Via Vetoio, Coppito

\noindent I-67010 L'Aquila (Italy)

\noindent e-mail ottaviani@vxscaq.aquila.infn.it

\end